\documentclass[10pt,twoside,reqno]{amsart}
\usepackage{amssymb}
\usepackage{multirow}
\calclayout

\numberwithin{equation}{section}
\makeatletter

\renewcommand{\@secnumfont}{\bfseries}

\renewcommand{\section}{\@startsection{section}{1}%
  {0mm}{.7\linespacing\@plus\linespacing}{.5\linespacing}
  {\normalfont\bfseries\centering}}

\newcommand{\bibsection}{\@startsection{section}{1}%
  {0mm}{.7\linespacing\@plus\linespacing}{.5\linespacing}
  {\normalfont\scshape\centering}}

\renewcommand{\@biblabel}[1]{#1.}
\newtheorem{theorem}{Theorem}[section]

\begin{document}

\title{Some identities of symmetry for the degenerate $q$-Bernoulli polynomials under symmetry group of degree $n$ }

\author{D. V. Dolgy}
\address{Hanrimwon, Kwangwoon University, Seoul 139-701, Republic of Korea\\Institute of Natural Sciences, Far Eastern Federal University, Vladivostok 690950, Russia.
}
\email{d\_dol@mail.ru}

\author{Taekyun Kim}
\address{Department of Mathematics, College of Science, Tianjin Polytechnic University, Tianjin City, 300387, China \\ Department of Mathematics, Kwangwoon University, Seoul 139-701, S. Korea}
\email{tkkim@kw.ac.kr}

\author{Lee-Chae Jang}
\address{Graduate School of Education, Konkuk University, Seoul 143-701, Republic of Korea}
\email{Icjang@konkuk.ac.kr}

\author{Hyuck-In Kwon}
\address{Department of Mathematics, Kwangwoon University, Seoul 139-701, Republic
of Korea}
\email{sura@kw.ac.kr}

\begin{abstract}
Recently, Kim-Kim Introduced some interesting identities of symmetry for $q$-Bernoulli polynomials under symmetry group of degree $n$. In this paper, we study the degenerate $q$-Euler polynomials and derive some identities of symmetry for these polynomials arising from the $p$-adic $q$-integral on $\mathbb{Z}_p$.
\end{abstract}

\maketitle
\global\long\def\acl#1#2{\left\langle \left.#1\right|#2\right\rangle }

\global\long\def\acr#1#2{\left\langle #1\left|#2\right.\right\rangle }

\global\long\def\Li{\mathrm{Li}}

\global\long\def\Zp{\mathbb{Z}_{p}}
\markboth{\centerline{\scriptsize Some identities of symmetry for the degenerate $q$-Bernoulli polynomials}}{\centerline{\scriptsize Some identities of symmetry for the degenerate $q$-Bernoulli polynomials}}

\section{Introduction}

Let $p$ be a fixed odd prime number. Throughout this paper, $\mathbb{Z}_{p}$, $\mathbb{Q}_{p}$ and $\mathbb{C}_{p}$ will denote the ring of $p$-adic integers, the field of $p$-adic
rational numbers and the completion of the algebraic closure of $\mathbb{Q}_{p}$.
The $p$-adic norm $|\cdot|_p$ is normalized as $|p|_p = \frac{1}{p}$.
Let $q \in \mathbb{C}_p$ be an indeterminate such that $\left|1-q\right|_{p}<p^{-\frac{1}{p-1}}$. The $q$-analogue of the number $x$ is defined by $\left[x\right]_{q}=\frac{1-q^{x}}{1-q}$. Let $f(x)$ be Uniformly differentiable function on $\mathbb{Z}_p$. The $p$-adic $q$-integral on $\mathbb{Z}_p$ is defined by Kim to be
\begin{equation}\begin{split}\label{01}
I_{q}(f) &= \int_{\mathbb{Z}_p} f(x) d\mu_{q}(x) = \lim_{N \rightarrow \infty} \sum_{x=0}^{p^N-1} f(x) \mu_{q} (x+p^N \mathbb{Z}_p) \\
&= \lim_{N \rightarrow \infty} \frac{1}{[p^N]_q} \sum_{x=0}^{p^N-1} f(x) q^x, \quad (\text{see}\,\,[13]).
\end{split}\end{equation}
In [1], L. Carlitz considered $q$-analogue of Bernoulli numbers which are given by recurrence relation to be
\begin{equation}\begin{split}\label{02}
\beta_{0,q}=1, \quad q(q\beta_q +1)^n - \beta_{n,q} = \begin{cases}1, \quad \text{if}\,\,\, n =1 \\ 0, \quad \text{if} \,\,\, n > 1, \end{cases}
\end{split}\end{equation}
with the usual convention about replacing $\beta_q^n$ by $\beta_{n,q}$. He defined $q$-Bernoulli polynomials as 
\begin{equation}\begin{split}\label{03}
\beta_{n,q}(x) = \sum_{l=0}^n {n \choose l} q^{lx} \beta_{l,q} [x]_q^{n-l}, \quad (\text{see} \,\, [1,13]).
\end{split}\end{equation}
In [13], Kim proved that the Carlitz's $q$-Bernoulli polynomials are represented as the $p$-adic $q$-integral on $\Zp$ which are given by 
\begin{equation}\begin{split}\label{04}
\int_{\mathbb{Z}_p}[x+y]_q^n    d\mu_{q} (y) = \beta_{n,q}(x), \quad (n \geq 0).
\end{split}\end{equation}
When $x=0$, $\beta_{n,q} = \beta_{n,q}(0)$ are the Carlitz $q$-Bernoulli numbers. 

In [2], L. Carlitz also introduced the degenerate Bernoulli polynomials which are given by the generating function to be 
\begin{equation}\begin{split}\label{05}
\frac{t}{(1+\lambda t)^{\frac{1}{\lambda}} -1} (1+\lambda t)^{\frac{x}{\lambda}} = \sum_{n=0}^\infty B_{n,\lambda}^* (x) \frac{t^n}{n!}.
\end{split}\end{equation}
Note that $\lim_{} B_{n,\lambda}^* (x) = B_{n}^* (x)$, where $B_n(x)$ are ordinary Bernoulli polynomials (see [1-10]). When $x=0$, $B_{n,\lambda}^* = B_{n,\lambda}^* (0)$ are called the degenerate Bernoulli numbers.
Recently, Kim-Kim introduced (fully) degenerate Bernoulli polynomials which are derived from
the $p$-adic invariant integral on $\Zp$ as follows:
\begin{equation}\begin{split}\label{06}
\int_{\mathbb{Z}_p}  (1+\lambda t)^{ \frac{x+y}{\lambda}}  d\mu_{q} (x) = \frac{\log(1+\lambda t)^{\frac{1}{\lambda}}}{(1+\lambda t)^{\frac{1}{\lambda}}-1} (1+\lambda t)^{\frac{x}{\lambda}}\quad (\textnormal{see} \,\, [7]),
\end{split}\end{equation}
where $\lambda, t \in \mathbb{C}_p$ with $|\lambda t|_p < p^{-\frac{1}{p-1}}$, and
\begin{equation*}\begin{split}
\lim_{q \rightarrow 1} \int_{\mathbb{Z}_p}  f(x)  d\mu_{q} (x) = \int_{\mathbb{Z}_p}  f(x)  d\mu_{1} (x) = \lim_{N \rightarrow \infty} \frac{1}{p^N} \sum_{x=0}^{p^N-1} f(x).
\end{split}\end{equation*}
The (fully) degenerate Bernoulli polynomials are defined by the generating function to be

\begin{equation}\begin{split}\label{07}
\frac{\log(1+\lambda t)}{\lambda(1+\lambda t)^{\frac{1}{\lambda}}-\lambda } (1+\lambda t)^{\frac{x}{\lambda}} = \sum_{n=0}^\infty B_{n,\lambda} (x) \frac{t^n}{n!}\quad (\textnormal{see} \,\, [7]).
\end{split}\end{equation}
Note that Kim's degenerate Bernoulli polynomials are slightly different from the Carlitz's degenerate Bernoulli polynomials. 

From \eqref{06} and \eqref{07}, we note that
\begin{equation}\begin{split}\label{08}
\lambda \int_{\mathbb{Z}_p} \left( \frac{x+y}{\lambda} \right)_n   d\mu_{1} (x) = B_{n,\lambda (x)}\quad (n \geq 0),
\end{split}\end{equation}
where $(x)_0 =1$, $(x)_n = x(x-1)\cdots (x-n+1), \,\, (n \geq 1)$. 

In [16], Kim considered degenerate $q$-Bernoulli polynomials which are given by the generating function to be 

\begin{equation}\begin{split}\label{09}
\int_{\mathbb{Z}_p} (1+ \lambda t)^{\frac{1}{\lambda}[x+y]_q}   d\mu_{q} (y) = \sum_{n=0}^\infty \beta_{n,\lambda,q} (x) \frac{t^n}{n!}.
\end{split}\end{equation}
When $x=0$, $\beta_{n,\lambda,q} = \beta_{n,\lambda,q}(0)$ are called (fully) degenerate $q$-Bernoulli numbers. Note that $\lim_{\lambda \rightarrow 0} \beta_{n,\lambda,q} (x) = \beta_{n,q}(x)$, $(n \geq 0)$.

In this paper, we give some identities of symmetry for the degenerate $q$-Bernoulli polynomials under symmetry group of degree $n$ arising from the $p$-adic $q$-integral on $\Zp$.

\section{Identities of symmetry for the degenerate $q$-Bernoulli polynomials}

We assume that $\lambda, t \in \mathbb{C}_p$ with $|\lambda|_p \leq 1, \,\,|t|_p < p^{-\frac{1}{p-1}}$. In this section, let  $w_1,w_2,\cdots,w_n$ be   positive integers. For $N \in \mathbb{N}$, we have
\begin{equation}\begin{split}\label{10}
&\int_{\mathbb{Z}_p} (1+ \lambda t)^{\frac{1}{\lambda} \big[ w_1 w_2 \cdots w_{n-1}y + w_1 \cdots w_n x + w_n \sum_{j=1}^{n-1}\big(\prod_{\substack{i=1\\i \neq j}}^{n-1} w_i \big) k_j\big]_q}   d\mu_{q^{w_1 w_2 \cdots w_{n-1}}} (y)\\
=&\lim_{N \rightarrow \infty} \frac{1}{[w_n p^N]_{q^{w_1 \cdots w_{n-1}}}} \\
&\times \sum_{y=0}^{w_n p^N -1} (1+ \lambda t)^{\frac{1}{\lambda} \big[ w_1 w_2 \cdots w_{n-1}y + w_1 \cdots w_n x + w_n \sum_{j=1}^{n-1}\big(\prod_{\substack{i=1\\i \neq j}}^{n-1} w_i \big) k_j\big]_q}  q^{w_1 w_2 \cdots w_{n-1}y}\\
=& \lim_{N \rightarrow \infty} \frac{1}{[w_n p^N]_{q^{w_1 \cdots w_{n-1}}}} \sum_{k_n=0}^{w_n-1} \sum_{y=0}^{p^N-1} q^{w_1 w_2 \cdots w_{n-1} (k_n + w_n y)}\\
&\times  (1+ \lambda t)^{\frac{1}{\lambda} \big[ \big( \sum_{j=1}^{n-1}w_j \big)(k_n +w_n y) + \sum_{j=1}^{n}w_j x + w_n \sum_{j=1}^{n-1}\big(\prod_{\substack{i=1\\i \neq j}}^{n-1} w_i \big) k_j\big]_q}. 
\end{split}\end{equation}
From \eqref{10}, we note that
\begin{equation}\begin{split}\label{11}
&\frac{1}{[w_1 \cdots w_{n-1}]_q}\prod_{l=1}^{n-1} \sum_{k_l=0}^{w_l-1} q^{w_n  \sum_{j=1}^{n-1}\big(\prod_{\substack{i=1\\i \neq j}}^{n-1} w_i \big)k_j}\\
&\times \int_{\mathbb{Z}_p} (1+ \lambda t)^{\frac{1}{\lambda} \big[ w_1 w_2 \cdots w_{n-1}y + w_1 \cdots w_n x + w_n \sum_{j=1}^{n-1}\big(\prod_{\substack{i=1\\i \neq j}}^{n-1} w_i \big) k_j\big]_q}   d\mu_{q^{w_1 w_2 \cdots w_{n-1}}} (y)\\
=& \lim_{N \rightarrow \infty} \frac{1}{[w_1 \cdots w_n p^N]_q}\prod_{l=1}^{n-1} \sum_{k_l=0}^{w_l-1} \sum_{k_n=0}^{w_n-1} \sum_{y=0}^{p^N-1} q^{w_1 w_2 \cdots w_{n-1} (k_n + w_n y) + \sum_{j=1}^{n-1}\big(\prod_{\substack{i=1\\i \neq j}}^{n-1} w_i \big) k_j w_n }\\
&\times (1+ \lambda t)^{\frac{1}{\lambda} \big[ \big( \sum_{j=1}^{n-1}w_j \big)(k_n +w_n y) + \sum_{j=1}^{n}w_j x + w_n \sum_{j=1}^{n-1}\big(\prod_{\substack{i=1\\i \neq j}}^{n-1} w_i \big) k_j\big]_q}. 
\end{split}\end{equation}
It is easy to show that \eqref{11} is invariant under any permutation in the symmetry group of degree $n$. Therefore, by \eqref{11}, we obtain the following theorem.\\\\

\begin{theorem}
Let $w_1, w_2, \cdots, w_n$ be positive integers . Then, the following expressions
\begin{equation*}\begin{split}
&\frac{1}{[w_{\sigma(1)} \cdots w_{\sigma(n-1)}]_q}\prod_{l=1}^{n-1} \sum_{k_l=0}^{w_{\sigma(l)}-1} q^{w_{\sigma(n)}  \sum_{j=1}^{n-1}\big(\prod_{\substack{i=1\\i \neq j}}^{n-1} w_{\sigma(i)} \big)k_j}\\
&\int_{\mathbb{Z}_p} (1+ \lambda t)^{\frac{1}{\lambda} \big[ w_{\sigma(1)} w_{\sigma(2)} \cdots w_{\sigma(n-1)} y + \sum_{j=1}^n w_j x + w_{\sigma(n)} \sum_{j=1}^{n-1}\big(\prod_{\substack{i=1\\i \neq j}}^{n-1} w_i \big) k_j\big]_q}   d\mu_{q^{w_{\sigma(1)} w_{\sigma(2)} \cdots w_{\sigma(n-1)}}} (y)
\end{split}\end{equation*}
are the same for any permutation $\sigma$ in the symmetry group of order $n$.
\end{theorem}
It is not difficult to show that
\begin{equation}\begin{split}\label{12}
&\bigg[ w_1 w_2 \cdots w_{n-1} y + w_1 w_2 \cdots w_n x + w_n \sum_{j=1}^{n-1}\big(\prod_{\substack{i=1\\i \neq j}}^{n-1} w_i \big) k_j \bigg]_q\\
&=[w_1 w_2 \cdots w_{n-1}]_q \Big[y + w_n x + \frac{w_n}{w_1} k_1+ \cdots + 
\frac{w_n}{w_{n-1}} k_{n-1} \Big]_{q^{w_1 w_2 \cdots w_{n-1}}}.
\end{split}\end{equation}
From \eqref{12}, we note that
\begin{equation}\begin{split}\label{13}
&\int_{\mathbb{Z}_p} (1+ \lambda t)^{\frac{1}{\lambda} \big[ w_1 \cdots w_{n-1}y + w_1 \cdots w_n x + w_n \sum_{j=1}^{n-1}\big(\prod_{\substack{i=1\\i \neq j}}^{n-1} w_i \big) k_j\big]_q}   d\mu_{q^{w_1 w_2 \cdots w_{n-1}}} (y)\\
=& \int_{\mathbb{Z}_p} (1+ \lambda t)^{\frac{[w_1 \cdots w_{n-1}]_q}{\lambda} \Big[y + w_n x + \frac{w_n}{w_1} k_1+ \cdots + \frac{w_n}{w_{n-1}} k_{n-1} \Big]_{q^{w_1 \cdots w_{n-1}}}}    d\mu_{q^{w_1 \cdots w_{n-1}}} (y)\\
=&  \int_{\mathbb{Z}_p} \left(1+ \frac{\lambda}{[w_1\cdots w_{n-1}]_q} [w_1 \cdots w_{n-1}]_q t \right)^{\frac{[w_1 \cdots w_{n-1}]_q}{\lambda} \Big[y + w_n x + \frac{w_n}{w_1} k_1+ \cdots + \frac{w_n}{w_{n-1}} k_{n-1} \Big]_{q^{w_1 \cdots w_{n-1}}}}  \\&\times   d\mu_{q^{w_1 \cdots w_{n-1}}} (y)\\
=& \sum_{m=0}^\infty [w_1 \cdots w_{n-1} ]_q^m \beta_{m, \frac{\lambda}{[w_1\cdots w_{n-1}]_q}, q^{w_1 \cdots w_{n-1}}} \Big(w_n x + \frac{w_n}{w_1} k_1+ \cdots + 
\frac{w_n}{w_{n-1}} k_{n-1} \Big) \frac{t^n}{n!}.
\end{split}\end{equation}
Therefore, by Theorem 1 and \eqref{13}, we obtain the following theorem.
\begin{theorem}
For $m \geq 0$, $w_1, w_2, \cdots, w_n \in \mathbb{N}$, the following expressions
\begin{equation*}\begin{split}
&[w_{\sigma(1)} \cdots w_{\sigma(n-1)} ]_q^{m-1} \prod_{l=1}^{n-1} \sum_{k_l=0}^{w_{\sigma(l)}-1} q^{\sum_{j=1}^{n-1}\big(\prod_{\substack{i=1\\i \neq j}}^{n-1} w_{\sigma(i)} \big) k_j w_{\sigma(n)}} \\
&\times \beta_{m, \frac{\lambda}{[w_{\sigma(1)}\cdots w_{\sigma(n-1)}]_q}, q^{w_{\sigma(1)} \cdots w_{\sigma(n-1)}}} \Big(w_{\sigma(n)} x + \frac{w_{\sigma(n)}}{w_{\sigma(1)}} k_1+ \cdots + \frac{w_{\sigma(n)}}{w_{\sigma(n-1)}} k_{n-1} \Big)
\end{split}\end{equation*}
are the same for any permutation $\sigma$ in the symmetry group of order $n$.
\end{theorem}
From \eqref{09}, we note that
\begin{equation}\begin{split}\label{14}
\sum_{n=0}^\infty \beta_{n,\lambda,q}(x) \frac{t^n}{n!} =&
\int_{\mathbb{Z}_p} (1+\lambda t)^{ \frac{1}{\lambda}[x+y]_q}   d\mu_{q} (y)\\
=& \sum_{n=0}^\infty \int_{\mathbb{Z}_p} \left( \frac{\frac{[x+y]_q}{\lambda}}{n} \right)  
 d\mu_{q} (x) \lambda^n t^n\\
=& \sum_{n=0}^\infty \lambda^n \int_{\mathbb{Z}_p} \left( \frac{[x+y]_q}{\lambda} \right)_n   d\mu_{q} (x) \frac{t^n}{n!}.
\end{split}\end{equation}
By comparing the coefficients on the both sides of \eqref{14}, we get
\begin{equation}\begin{split}\label{15}
\beta_{n,\lambda,q} =& \lambda^n \int_{\mathbb{Z}_p} \left( \frac{[x+y]_q}{\lambda} \right)_n   d\mu_{q} (x)\\
=&\lambda^n \sum_{m=0}^n S_1 (n,m) \lambda^{-m} \int_{\mathbb{Z}_p} [x+y]_q^m   d\mu_{q} (y)\\
=& \sum_{m=0}^n S_1(n,m) \lambda^{n-m} \beta_{m,q}(x).
\end{split}\end{equation}
where $\beta_{m,q}(x)$ are called Carlitz's $q$-Bernoulli polynomials.

Now, we observe that

\begin{equation}\begin{split}\label{16}
&\Big[y + w_n x + w_n \sum_{j=1}^{n-1} \frac{k_j}{w_j} \Big]_{q^{w_1 \cdots w_{n-1}}}
\\
&= \frac{[w_n]_q}{[w_1 \cdots w_{n-1}]_q} \bigg[ \sum_{j=1}^{n-1}\big(\prod_{\substack{i=1\\i \neq j}}^{n-1} w_i \big) k_j \bigg]_{q^{w_n}}+q^{w_n \sum_{j=1}^{n-1}\big(\prod_{\substack{i=1\\i \neq j}}^{n-1} w_i \big) k_j}[y+w_n x]_{q^{w_1 \cdots w_{n-1}}}.
\end{split}\end{equation}
By \eqref{15}, we get
\begin{equation}\begin{split}\label{17}
&\beta_{m, \frac{\lambda}{[w_1\cdots w_{n-1}]_q}, q^{w_1 \cdots w_{n-1}}} \Big(w_n x + \frac{w_n}{w_1} k_1+ \cdots + \frac{w_n}{w_{n-1}} k_{n-1} \Big)\\
=& \left( \frac{\lambda}{[w_1 \cdots w_{n-1} ]_q} \right)^m \int_{\mathbb{Z}_p} \left( \Big(\frac{\lambda}{[w_1 \cdots w_{n-1}]_q} \Big)^{-1} \Big[y + w_n x + w_n \sum_{j=1}^{n-1} \frac{k_j}{w_j} \Big]_q^l \right)   d\mu_{q^{w_1 \cdots w_{n-1}} } (y) \\
=&\left( \frac{\lambda}{[w_1 \cdots w_{n-1} ]_q} \right)^m \sum_{l=0}^m S_1 (m,l) [w_1 \cdots w_{n-1}]_q^l \lambda^{-l} \\
&\quad \times \int_{\mathbb{Z}_p} \Big[y + w_n x + w_n \sum_{j=1}^{n-1} \frac{k_j}{w_j} \Big]_q^l  d\mu_{q^{w_1 \cdots w_{n-1}} } (y).
\end{split}\end{equation}
From \eqref{16}, we can derive the following equation:
\begin{equation}\begin{split}\label{18}
&\int_{\mathbb{Z}_p} \Big[y + w_n x + w_n \sum_{j=1}^{n-1} \frac{k_j}{w_j} \Big]_q^l  d\mu_{q^{w_1 \cdots w_{n-1}} } (y) \\
=&\sum_{s=0}^l {l \choose s} \left( \frac{[w_n]_q}{[w_1 \cdots w_{n-1}]_q} \right)^{l-s} \bigg[ \sum_{j=1}^{n-1}\big(\prod_{\substack{i=1\\i \neq j}}^{n-1} w_i \big) k_j \bigg]_{q^{w_n}}^{l-s} q^{w_n s\sum_{j=1}^{n-1}\big(\prod_{\substack{i=1\\i \neq j}}^{n-1} w_i \big) k_j}\\
&\times  \int_{\mathbb{Z}_p} [y+w_n x]_{q^{w_1 \cdots w_{n-1}}}^s  d\mu_{q^{w_1 \cdots w_{n-1}} } (y) \\
=& \sum_{s=0}^l {l \choose s} \left( \frac{[w_n]_q}{[w_1 \cdots w_{n-1}]_q} \right)^{l-s} \bigg[ \sum_{j=1}^{n-1}\big(\prod_{\substack{i=1\\i \neq j}}^{n-1} w_i \big) k_j \bigg]_{q^{w_n}}^{l-s} q^{w_n s\sum_{j=1}^{n-1}\big(\prod_{\substack{i=1\\i \neq j}}^{n-1} w_i \big) k_j} \\ &\times \beta_{s, q^{w_1 \cdots w_{n-1}}} \big(w_n x \big).
\end{split}\end{equation}
By \eqref{17} and \eqref{18}, we get
\begin{equation}\begin{split}\label{19}
&\beta_{m, \frac{\lambda}{[w_1\cdots w_{n-1}]_q}, q^{w_1 \cdots w_{n-1}}} \Big(w_n x + \frac{w_n}{w_1} k_1+ \cdots + \frac{w_n}{w_{n-1}} k_{n-1} \Big)\\
=& \sum_{p=0}^m \sum_{s=0}^p {p \choose s} S_1(m,p) \lambda^{m-p} [w_1 \cdots w_{n-1}]_q^{s-m} [w_n]_q^{p-s} \bigg[ \sum_{j=1}^{n-1}\big(\prod_{\substack{i=1\\i \neq j}}^{n-1} w_i \big) k_j \bigg]_{q^{w_n}}^{p-s}\\
& \times q^{w_n s\sum_{j=1}^{n-1}\big(\prod_{\substack{i=1\\i \neq j}}^{n-1} w_i \big) k_j}\beta_{s, q^{w_1 \cdots w_{n-1}}} \big(w_n x \big).
\end{split}\end{equation}
From \eqref{19}, we note that
\begin{equation}\begin{split}\label{20}
&[w_1 \cdots w_{n-1}]_q^{m-1} \prod_{l=1}^{n-1} \sum_{k_l=0}^{w_l-1} q^{w_n  \sum_{j=1}^{n-1}\big(\prod_{\substack{i=1\\i \neq j}}^{n-1} w_i \big)k_j} \\
&\times \beta_{m, \frac{\lambda}{[w_1\cdots w_{n-1}]_q}, q^{w_1 \cdots w_{n-1}}} \Big(w_n x + w_n \sum_{j=1}^{n-1} \frac{k_j}{w_j} \Big)\\
=& \prod_{l=1}^{n-1} \sum_{k_l=0}^{w_l-1}\sum_{p=0}^m \sum_{s=0}^p {p \choose s} S_1(m,p) \lambda^{m-p} [w_1 \cdots w_{n-1}]_q^{s-1} [w_n]_q^{p-s}
\bigg[ \sum_{j=1}^{n-1}\big(\prod_{\substack{i=1\\i \neq j}}^{n-1} w_i \big) k_j \bigg]_{q^{w_n}}^{p-s}\\
&\times q^{(s+1) w_n \sum_{j=1}^{n-1}\big(\prod_{\substack{i=1\\i \neq j}}^{n-1} w_i \big) k_j}  \beta_{s, q^{w_1 \cdots w_{n-1}}} \big(w_n x \big)\\
=& \sum_{p=0}^m \sum_{s=0}^p {p \choose s} S_1(m,p) \lambda^{m-p} [w_1 \cdots w_{n-1}]_q^{s-1} [w_n]_q^{p-s} \beta_{s, q^{w_1 \cdots w_{n-1}}} \big(w_n x \big)\\
&\times \prod_{l=1}^{n-1} \sum_{k_l=0}^{w_l-1}q^{(s+1) w_n \sum_{j=1}^{n-1}\big(\prod_{\substack{i=1\\i \neq j}}^{n-1} w_i \big) k_j}
\bigg[ \sum_{j=1}^{n-1}\big(\prod_{\substack{i=1\\i \neq j}}^{n-1} w_i \big) k_j \bigg]_{q^{w_n}}^{p-s}\\
=& \sum_{p=0}^m \sum_{s=0}^p {p \choose s} S_1(m,p) \lambda^{m-p} [w_1 \cdots w_{n-1}]_q^{s-1} [w_n]_q^{p-s} \beta_{s, q^{w_1 \cdots w_{n-1}}} \big(w_n x \big)\\
&\times K_{n,q^{w_n}}(w_1, \cdots w_{n-1}| p-s,s),
\end{split}\end{equation}
where
\begin{equation}\begin{split}\label{21}
 K_{n,q}(w_1, \cdots w_{n-1}|i,t) = \prod_{l=1}^{n-1} \sum_{k_l=0}^{w_l-1}q^{(t+1) \sum_{j=1}^{n-1}\big(\prod_{\substack{i=1\\i \neq j}}^{n-1} w_i \big) k_j}
\bigg[ \sum_{j=1}^{n-1}\big(\prod_{\substack{i=1\\i \neq j}}^{n-1} w_i \big) k_j \bigg]_q^i .
\end{split}\end{equation}
Therefore, by \eqref{20} and \eqref{21}, we obtain the following theorem.

\begin{theorem}
Let $m \geq 0$ and $w_1, w_2, \cdots, w_n \mathbb{N}$, Then the following expressions
\begin{equation*}\begin{split}
&\sum_{p=0}^m \sum_{s=0}^p {p \choose s} S_1(m,p) \lambda^{m-p} [w_{\sigma(1)} \cdots w_{\sigma(n-1)}]_q^{s-1} [w_{\sigma(n)}]_q^{p-s} \beta_{s, q^{w_{\sigma(1)} \cdots w_{\sigma(n-1)}}} \big(w_{\sigma(n)} x \big)\\
&\times K_{n,q^{w_{\sigma(n)}}}(w_{\sigma(1)}, \cdots w_{\sigma(n-1)}| p-s,s)
\end{split}\end{equation*}
are the same for any permutation $\sigma$ in the symmetry group of order $n$.
\end{theorem}
Note that some identities of Bernoulli and Euler polynomials are studied by several authors (see [1-19]).

\noindent
{\bf{Acknowledgements}}\\\\
 This paper is supported by grant NO 14-11-00022 of Russian Scientific Fund.

\bibliographystyle{amsplain}

\begin{thebibliography}{10}

\bibitem{key-1} L. Carlitz, \textit{$q$-Bernoulli and Eulerian numbers
}, Trans. Amer. Math. Soc., {\bf{76}}  (1954), 332--350.


\bibitem{key-2} L. Carlitz, \textit{Degenerate Stirling, Bernoulli and Eulerian numbers}, Utilitas Math., {\bf{15}} (1979), 51--88.


\bibitem{key-3} Y. He,  \textit{Symmetric identities for Carlitz's $q$-Bernoulli numbers and polynomials}, Adv. Difference Equ., {\bf{2013}} 2013:246, 10 pp. 

\bibitem{key-4} D. S. Kim, T. Kim,
\textit{Some identities of symmetry for $q$-Bernoulli polynomials under symmetric group of degree $n$}, Ars Comb., {\bf{126}} (2016), 435--441.

\bibitem{key-5} D. S. Kim, N. Lee, J. Na, K. H. Park,
\textit{Abundant symmetry for higher-order Bernoulli polynomials (I)
}, Adv. Stud. Contemp. Math., {\bf{23}} (2013), no. 3, 461--482. 

\bibitem{key-6} D. S. Kim, N. Lee, J. Na, K. H. Park,
\textit{Identities of symmetry for higher-order Euler polynomials in three variables (I)
}, Adv. Stud Contemp. Math., {\bf{22}} (2012), no. 1, 51--74.

\bibitem{key-7} T. Kim, D. S. Kim, J.-J. Seo, \textit{Fully degenerate poly-Bernoulli numbers and polynomials}, Open Math., {\bf{2016}}; 14: 545-556. 

\bibitem{key-8} T. Kim, D. V. Dolgy, J. J. Seo,
\textit{Identities of symmetry for degenerate $q$-Euler polynomials}, Adv. Stud. Contemp. Math., {\bf{25}} (2015), 577--582.

\bibitem{key-9} T. Kim, D. V. Dolgy, D. S. Kim,
\textit{Symmetric identities for degenerate generalized Bernoulli polynomials
}, J. Nonlinear Sci. Appl., {\bf{9}} (2016), no. 2, 677--683.

\bibitem{key-10} T. Kim,
\textit{Symmetric identities of degenerate Bernoulli polynomials
}, Proc. Jangjeon Math. Soc., {\bf{18}} (2015), no. 4, 593-599.

\bibitem{key-11} T. Kim, D. S. Kim, H.-I. Kwon, D. V. Dolgy,
\textit{Some identities of $q$-Euler polynomials under the symmetric group of degree $n$
}, J. Nonlinear Sci. Appl., {\bf{9}} (2016), no. 3, 1077--1082.

\bibitem{key-12} T.Kim  \textit{An identity of the symmetry for the Frobenius-Euler polynomials associated with the fermionic $p$-adic invariant $q$-integrals on $\bold Z\sb p$}, Rocky Mountain J. Math., {\bf{41}} (2011), no. 1, 239--247.
 
\bibitem{key-13} T. Kim, \textit{$q$-Volkenborn integration}, Russ. J. Math. Phys., {\bf{9}} (2002), no. 3, 288-299. 

\bibitem{key-14} T. Kim, H.-I. Kwon, J.-J. Seo,
\textit{Identities of symmetry for degenerate $q$-Bernoulli polynomials}, Proc. Jangjeon Math. Soc., {\bf{18}} (2015), no. 4, 495--499.


\bibitem{key-15} T. Kim,
\textit{Some identities of the $q$-Euler polynomials of higher-order and $q$-Stirling numbers by the fermionic $p$-adic integral on $\mathbb{Z}_p$}, Russ. J. Math. Phys., {\bf{16}} (2009), no. 4, 484--491.


\bibitem{key-16} T. Kim, \textit{On degenerate $q$-Bernoulli polynomials}, Bull. Korean Math. Soc. {\bf{53}} (2016), no. 4, 1149-1156.

\bibitem{key-17} Y.-H. Kim, K.-W. Hwang,
\textit{Symmetry of power sum and twisted Bernoulli polynomials},
Adv. Stud. Contemp. Math., {\bf{18}} (2009), no. 2, 127--133.

\bibitem{key-18} H. I. Kwon, T. Kim, J. J. Seo,
\textit{Some identities of symmetry for modified degenerate Frobenius-Euler polynomials},
Adv. Stud. Contemp. Math., {\bf{26}} (2016), 299--305.

\bibitem{key-19} E.-J. Moon, S.-H. Rim, J.-H. Jin, S.-J. Lee,
\textit{On the symmetric properties of higher-order twisted $q$-Euler numbers polynomials},
Adv. Difference Equ., {\bf{2010}} Art. ID 765259, 8 pages.

\end{thebibliography}
\providecommand{\bysame}{\leavevmode\hbox to3em{\hrulefill}\thinspace}
\providecommand{\MR}{\relax\ifhmode\unskip\space\fi MR }
\providecommand{\MRhref}[2]{%
  \href{http://www.ams.org/mathscinet-getitem?mr=#1}{#2}
}
\providecommand{\href}[2]{#2}

\end{document}